\documentclass[12pt,reqno]{amsart}
\usepackage{amsmath,amssymb,amsfonts,amscd,latexsym,amsthm,mathrsfs,verbatim,comment,cite}
\usepackage[unicode]{hyperref}
\newtheorem{theorem}{Theorem}

\theoremstyle{remark}

\newtheorem*{remarks}{Remarks}
\begin{document}

\title{A remarkable basic hypergeometric identity}

\author{Christian Krattenthaler}
\address{Fakult\"at f\"ur Mathematik, Universit\"at Wien, Oskar-Morgenstern-Platz 1, A-1090 Vienna, Austria}
\urladdr{https://www.mat.univie.ac.at/~kratt/}

\author{Wadim Zudilin}
\address{Department of Mathematics, IMAPP, Radboud University, PO Box 9010, 6500~GL Nijmegen, The Netherlands}
\urladdr{https://www.math.ru.nl/~wzudilin/}

\date{5 June 2024}

\dedicatory{To George Andrews and Bruce Berndt, with admiration,\\on their jointly complex $|2+9\sqrt{-1}|^2$th birthday}



\begin{abstract}
We give a closed form for \emph{quotients} of truncated basic hypergeometric series where the base~$q$ is evaluated at roots of unity.
\end{abstract}

\thanks{The first author was partially supported by the Austrian Science Fund FWF, grant 10.55776/F50.}

\maketitle

We set $(a;q)_n=(1-a)(1-aq)\dotsb(1-aq^{n-1})$ to be the $q$-shifted factorial ($q$-Pochhammer symbol), with its multiple version
\[
(a_1,\dots,a_m;q)_k=\prod_{j=1}^m(a_j;q)_k.
\]

An investigation of supercongruences for truncated $_4F_3$ hypergeometric sums in \cite{Zu20} revealed an interesting structure for their $q$-counterparts at roots of unity.
In order to state the corresponding identity, we introduce the sum
\begin{equation}
\begin{gathered}
F_n(a,q;\ell_1,\ell_2)
=\sum_{k=0}^{n-1}f_k(a,q;\ell_1,\ell_2),
\\ 
\text{where}\;
f_k(a,q;\ell_1,\ell_2)=\frac{(q^{\ell_1}a,q^{1-\ell_1}a,q^{\ell_2}a,q^{1-\ell_2}a;q)_k}{(qa;q)_k^4}\,q^k,
\end{gathered}
\label{eq:F}
\end{equation}
and the product
\begin{equation}
G(a,q;\ell_1,\ell_2)
=\frac{\prod_{j=0}^{\ell_1-1}(a-q^j)\cdot\prod_{j=0}^{\ell_2-1}(a-q^j)}
{\prod_{j=0}^{\ell_1-1}(1-q^ja)\cdot\prod_{j=0}^{\ell_2-1}(1-q^ja)}.
\label{eq:G}
\end{equation}
Here $\ell_1,\ell_2$ are integers; they are assumed to be non-negative in~\eqref{eq:G}. If one (or both) of $\ell_1$ and $\ell_2$ are zero, empty products
have to be interpreted as~1.

\begin{theorem}
\label{ob1}
For $\zeta=\zeta_n$ a primitive $n$th root of unity and $\ell_1,\ell_2\ge0$ integers, define $F_n(a,\zeta)=F_n(a,\zeta;\ell_1,\ell_2)$.
Then
\begin{equation}
\frac{F_n(a,\zeta)}{F_n(1,\zeta)}
=\frac{n^2a^{n-1}}{(1+a+a^2+\dots+a^{n-1})^2}\,G(a,\zeta;\ell_1,\ell_2).
\label{eq:id}
\end{equation}
\end{theorem}

\begin{remarks}
(1) As a corollary, we obtain
\[
\frac{F_n(a,\zeta)F_n(1/a,\zeta)}{F_n(1,\zeta)^2}=\bigg(\frac{na^{(n-1)/2}}{1+a+a^2+\dots+a^{n-1}}\bigg)^4
\]
independent of the choice of the primitive $n$th root of unity;
this formula is a result stated without proof in \cite{Zu20}.

\medskip
(2)
We may extend the definition of product to cases where the upper bound
on the product index is less than the lower bound by
\begin{equation} \label{eq:prod}
\prod _{k=M} ^{N-1}\operatorname{Expr}(k)=\begin{cases} \hphantom{-}
\prod _{k=M} ^{N-1} \operatorname{Expr}(k)&N>M,\\
\hphantom{-}1&N=M\\
1\Big/\prod _{k=N} ^{M-1}\operatorname{Expr}(k)&N<M.\end{cases}
\end{equation}
With this convention, we may observe that
$$
G(a,q;-\ell_1,\ell_2)=G(a,q;\ell_1+1,\ell_2),
$$
and the analogous equation if the roles of $\ell_1$ and $\ell_2$ are
interchanged. Since, trivially, $F_n(a,q;\ell_1,\ell_2)$ satisfies
the same identities, we see that the restriction on~$\ell_1$ and~$\ell_2$
in Theorem~\ref{ob1} can be dropped when assuming the convention
in~\eqref{eq:prod}.
\end{remarks}

The identity in Theorem~\ref{ob1} is somewhat unusual when compared to other evaluations of basic hypergeometric series, when those are expressed in `closed form'\,---\,a product or a finite linear combination of products.
We were not able to produce a closed-form evaluation for the sum $F_n(a,\zeta)$ itself, only for its renormalisation $F_n(a,\zeta)/F_n(1,\zeta)$ which is not a sum but a \emph{quotient} of sums.
Identity~\eqref{eq:id} connects visually\,---\,though not explicitly\,---\,with known evaluations involving the cyclic quantum dilogarithm \cite{KMS93} (see also \cite[Appendix~A]{GZ21});
it also shares similarities with identities obtained in~\cite{SZ97}.

\medskip
We prove Theorem~\ref{ob1} by constructing explicit difference equations for the sums $F_n(a,\zeta;\ell_1,\ell_2)$ with respect to the $\ell$-parameters, separately for the cases $\ell_1\ne\ell_2$ and $\ell_1=\ell_2$, by demonstrating that the expressions $G(a,\zeta;\ell_1,\ell_2)F_n(1,\zeta;\ell_1,\ell_2)$ satisfy the same equations, and by finally verifying the proportionality of these two solutions for a few starting values (namely, for the case $\ell_1=1$ and $\ell_2$~arbitrary).

\medskip
In the formal identity
\begin{multline*}
(d-b)(1-aK)(1-cK) + (a-d)(1-bK)(1-cK) \\ + (b-c)(1-aK)(1-dK) + (c-a)(1-bK)(1-dK) = 0,
\end{multline*}
multiply both sides by
\[
\frac{(a,b,c,d;q)_k}{e\cdot(qe;q)_k^4}\,q^k
\]
and replace $K\to q^k$, $a\to q^{\ell_1}a$, $b\to q^{-\ell_1}a$, $c\to q^{\ell_2}a$, $d\to q^{-\ell_2}a$, $e\to a$ to get
\begin{align*}
&
(q^{-\ell_2}-q^{-\ell_1})(1-q^{\ell_1}a)(1-q^{\ell_2}a) f_k(a,q;\ell_1+1,\ell_2+1)
\\ &\qquad
+ (q^{\ell_1}-q^{-\ell_2})(1-q^{-\ell_1}a)(1-q^{\ell_2}a) f_k(a,q;\ell_1,\ell_2+1)
\\ &\qquad
+ (q^{-\ell_1}-q^{\ell_2})(1-q^{\ell_1}a)(1-q^{-\ell_2}a) f_k(a,q;\ell_1+1,\ell_2)
\\ &\qquad
+ (q^{\ell_2}-q^{\ell_1})(1-q^{-\ell_1}a)(1-q^{-\ell_2}a) f_k(a,q;\ell_1,\ell_2)
= 0.
\end{align*}
Setting $a=1$ in the result, multiplying both sides by $G(a,q;\ell_1,\ell_2)$ and using
\[
G(a,q;\ell_1+1,\ell_2)
=\frac{a-q^{\ell_1}}{1-q^{\ell_1}a}\,G(a,q;\ell_1,\ell_2)
\]
and
\[
G(a,q;\ell_1,\ell_2+1)
=\frac{a-q^{\ell_2}}{1-q^{\ell_2}a}\,G(a,q;\ell_1,\ell_2)
\]
for the product in~\eqref{eq:G}, we find out that
\begin{align*}
&
(q^{-\ell_2}-q^{-\ell_1})(1-q^{\ell_1}a)(1-q^{\ell_2}a) g_k(a,q;\ell_1+1,\ell_2+1)
\\ &\qquad
+ (q^{\ell_1}-q^{-\ell_2})(1-q^{-\ell_1}a)(1-q^{\ell_2}a) g_k(a,q;\ell_1,\ell_2+1)
\\ &\qquad
+ (q^{-\ell_1}-q^{\ell_2})(1-q^{\ell_1}a)(1-q^{-\ell_2}a) g_k(a,q;\ell_1+1,\ell_2)
\\ &\qquad
+ (q^{\ell_2}-q^{\ell_1})(1-q^{-\ell_1}a)(1-q^{-\ell_2}a) g_k(a,q;\ell_1,\ell_2)
= 0,
\end{align*}
where
\[
g_k(a,q;\ell_1,\ell_2)
= G(a,q;\ell_1,\ell_2) f_k(1,q;\ell_1,\ell_2).
\]
In particular, both sums
\[
F_n(a,q;\ell_1,\ell_2) = \sum_{k=0}^{n-1} f_k(a,q;\ell_1,\ell_2)
\]
and
\[
\sum_{k=0}^{n-1} g_k(a,q;\ell_1,\ell_2) = G(a,q;\ell_1,\ell_2) F_n(a,q;\ell_1,\ell_2)
\]
satisfy the \emph{same} recursion
\begin{align}
&
(q^{-\ell_2}-q^{-\ell_1})(1-q^{\ell_1}a)(1-q^{\ell_2}a) F(a,q;\ell_1+1,\ell_2+1)
\nonumber\\ &\qquad
+ (q^{\ell_1}-q^{-\ell_2})(1-q^{-\ell_1}a)(1-q^{\ell_2}a) F(a,q;\ell_1,\ell_2+1)
\nonumber\\ &\qquad
+ (q^{-\ell_1}-q^{\ell_2})(1-q^{\ell_1}a)(1-q^{-\ell_2}a) F(a,q;\ell_1+1,\ell_2)
\nonumber\\ &\qquad
+ (q^{\ell_2}-q^{\ell_1})(1-q^{-\ell_1}a)(1-q^{-\ell_2}a) F(a,q;\ell_1,\ell_2)
= 0.
\label{eq4}
\end{align}
Note also that the second sum is shorter when $0<\ell_1,\ell_2<n$:
\begin{equation}
\sum_{k=0}^{n-1} g_k(a,q;\ell_1,\ell_2)
= G(a,q;\ell_1,\ell_2) \sum_{k=0}^{\min\{\ell_1,\ell_2\}-1} \frac{(q^{\ell_1},q^{1-\ell_1},q^{\ell_2},q^{1-\ell_2};q)_k}{(q;q)_k^4}\,q^k.
\label{eq:short}
\end{equation}
Furthermore observe that the difference equation in~\eqref{eq4} degenerates when $\ell_1=\ell_2$ (and therefore does not give access to comparison of the two solutions along the diagonal $\ell_1=\ell_2$).

\medskip
In the case $\ell_1=\ell_2$ the difference equation is three-term but less friendly.
For simplicity we write $\ell$ instead of $\ell_1=\ell_2$.
The announced difference equation corresponds to the difference operator
\begin{align*}
\mathcal{L}_q
&:=
- qL^2(1+L)(1+4L+L^2)(1-qL)^3(1-La)^2(1-qLa)^2 S^2
\\ &\;
+ (1-qL^2) (1 + 5(1+q)L + (5+16q+5q^2)L^2 + (1+q)(1-5q+q^2)L^3
\\ &\;\quad
- 2q(3+25q+3q^2)L^4 + q(1+q)(1-5q+q^2)L^5  + q^2(5+16q+5q^2)L^6
\\ &\;\quad
+ 5q^3(1+q)L^7 + q^4L^8) (1-La)^2(a-qL)^2 S
\\ &\;
- qL^2(1-L)^3(1+qL)(1+4qL+q^2L^2)(a-L)^2(a-qL)^2,
\end{align*}
where $S$ denotes the shift operator $\ell\mapsto\ell+1$, while $L=q^\ell$.
Explicitly, if we write $\tilde f_k(a,q;\ell)=s(a,q,q^\ell,q^ka) f_k(a,q;\ell,\ell)$ with
\begin{align*}
s(a,q,L,K)
&=\frac{q(1+L)(1+qL)(1-qL^2)(a-L)^2(a-qL)^2L^2(1-K)^4}{K(K-qL)^2(K-L)^2}
\\ &\;\times
(K(1+q^3L^6)+4K(1+q)(1+q^2L^4)L
\\ &\;\quad
-(4q^2-K-13qK-q^2K+4K^2)(1+qL^2)L^2
\\ &\;\quad
-2(q^3+7q(q+K^2)+K^2)L^3),
\end{align*}
then\footnote{We found this creative telescoping relation using the
  implementation of the $q$-analogue of Zeilberger's algorithm
  (cf.\ \cite{Z-Koorn,Ze90}) within the \texttt{QDifferenceEquations} package in~\textsl{Maple}.}
\[
\mathcal{L}_qf_k(a,q;\ell,\ell)=\tilde f_{k+1}(a,q;\ell)-\tilde f_k(a,q;\ell).
\]
Since $\tilde f_n(a,\zeta;\ell)=\tilde f_0(a,\zeta;\ell)$ for a primitive $n$th root of unity, we find by telescoping that
$\mathcal{L}_\zeta F_n(a,\zeta;\ell,\ell)=0$.
Because of the special shape of the difference operator $\mathcal{L}_q$ and of the certificate $s(a,q,L,K)$, similarly to the situation where $\ell_1\ne\ell_2$ we conclude that also
\[
\mathcal{L}_\zeta \big(G(a,\zeta;\ell,\ell) F_n(1,\zeta;\ell,\ell)\big)=0.
\]

\medskip
To prove the identity in Theorem~\ref{ob1}, we check it for $\ell_1=\ell_2=0$ and for $\ell_1=1$, $\ell_2$ arbitrary.
The cases $\ell_1=\ell_2=0$ and $\ell_1=\ell_2=1$ and the homogeneous recursion associated with the operator $\mathcal{L}_\zeta$ implies that \eqref{eq:id} holds when $\ell_1=\ell_2$ is arbitrary.
Then the other difference equation implies the validity of \eqref{eq:id} for $\ell_1=1$, for $\ell_1=\ell_2$, and the symmetry $\ell_1\leftrightarrow\ell_2$ of both sides implies the desired identity for general $\ell_1$ and $\ell_2$.

By definition, we have
\[
G(a,\zeta;0,0)
=G(a,\zeta;1,1)=1.
\]
At the same time, the equality $F_n(a,\zeta;0,0)=F_n(a,\zeta;1,1)$ follows from the definition~\eqref{eq:F}.
On the basis of the three-term recursion satisfied by $F_n(a,\zeta;\ell,\ell)$ and by $G(a,\zeta;\ell,\ell) F_n(1,\zeta;\ell,\ell)$ this implies the proportionality of the two for every $\ell\in\mathbb Z$.

When $\ell_1=1$, the sum $\sum_k g_k(a,\zeta;1,\ell_2)$ reduces to the single entry for $k=0$ (see~\eqref{eq:short}). Thus, this case requires
the verification of
\begin{align}
\sum_{k=0}^{n-1} \frac{(1-a)\,(\zeta^{\ell}a,\zeta^{1-\ell}a;\zeta)_k}{(1-\zeta^ka)\,(\zeta a;\zeta)_k^2}\,\zeta^k
&=\frac{n^2a^{n-1}}{(1+a+a^2+\dots+a^{n-1})^2}\,G(a,\zeta;1,\ell)
\nonumber\\
&=\frac{n^2a^{n-1}}{(1+a+a^2+\dots+a^{n-1})^2}\,\frac{\prod_{j=1}^{\ell-1}(a-\zeta^j)}{\prod_{j=1}^{\ell-1}(1-\zeta^ja)}.
\label{eq5}
\end{align}
We now take
\[
h_k(a,q;\ell)=\frac{(1-a)\,(q^{\ell}a,q^{1-\ell}a;q)_k}{(1-q^ka)\,(qa;q)_k^2}\,q^k
\]
and
\begin{align*}
\tilde h_k(a,q;\ell)
&=\frac{(1-q^ka)^3(1+q^\ell)(a-q^\ell)q^{\ell-k}}{a\,(1-q^\ell)^2(q^\ell-q^ka)}\,h_k(a,q;\ell)
\end{align*}
with the motive that
\[
(1-q^\ell a)h_k(a,q;\ell+1)
-(a-q^\ell)h_k(a,q;\ell)
=\tilde h_{k+1}(a,q;\ell)-\tilde h_k(a,q;\ell).
\]
Substituting $q=\zeta$ and summing over $k$ from $0$ to $n-1$ we obtain
\[
(1-\zeta^\ell a)H(a,\zeta;\ell+1)
-(a-\zeta^\ell)H(a,\zeta;\ell)
=\tilde h_n(a,\zeta;\ell)-\tilde h_0(a,\zeta;\ell)=0,
\]
where
\[
H(a,\zeta;\ell)=\sum_{k=0}^{n-1}h_k(a,\zeta;\ell)
=\sum_{k=0}^{n-1} \frac{(1-a)\,(\zeta^{\ell}a,\zeta^{1-\ell}a;\zeta)_k}{(1-\zeta^ka)\,(\zeta a;\zeta)_k^2}\,\zeta^k.
\]
This shows that verification in \eqref{eq5} reduces to the one where $\ell=1$, that is, to
\[
\sum_{k=0}^{n-1} \frac{1}{(1-\zeta^ka)^2}\,\zeta^k
=\frac{n^2a^{n-1}}{(1-a^n)^2}
=\frac{n^2a^{n-1}}{\prod_{j=0}^{n-1}(1-\zeta^ja)^2};
\]
the latter sum on the left-hand side is the partial fraction decomposition of the right-hand side.

\end{document}